\documentclass[twocolumn,letterpaper]{IEEEAerospaceCLS}  

\usepackage[]{graphicx,float,latexsym,amssymb,amsfonts,amsmath,amstext,times,tabularx,multirow,stfloats} 

\usepackage[]{graphicx}    
\newcommand{\ignore}[1]{}  

\begin{document}
\title{Determining Follower Aircraft's Optimal Trajectory in Relation to a Dynamic Formation Ring}

\author{%
Carl A. Gotwald \& Michael D. Zollars\\ 
Department of Aeronautics and Astronautics\\
Air Force Institute of Technology\\
Wright-Patterson AFB, OH 45433\\
carl.gotwald@afit.edu; michael.zollars@afit.edu
\and 
Isaac E. Weintraub\\
Controls Science System Center\\
Wright-Patterson AFB, OH 45433\\
isaac.weintraub.1@us.af.mil
\thanks{{U.S. Government work not protected by U.S. copyright}}         
}

\maketitle

\thispagestyle{plain}
\pagestyle{plain}

\maketitle

\thispagestyle{plain}
\pagestyle{plain}

\begin{abstract}
The specific objective of this paper is to develop a tool that calculates the optimal trajectory of the follower aircraft as it completes a formation rejoin, and then maintains the formation position, defined as a ring of points, until a fixed final time. The tool is designed to produce optimal trajectories for a variety of initial conditions and leader trajectories. Triple integrator dynamics are used to model the follower aircraft in three dimensions. Control is applied directly to the rate of acceleration. Both the follower’s and leader’s velocities and accelerations are bounded, as dictated by the aircraft’s performance envelope. Lastly, a path constraint is used to ensure the follower avoids the leader's jet wash region. This optimal control problem is solved through numerical analysis using the direct orthogonal collocation solver GPOPS-II. Two leader trajectories are investigated, including a descending spiral and continuous vertical loops. Additionally, a study of the effect of various initial guesses is performed. All trajectories displayed a direct capture of the formation position, however changes in solver initial conditions demonstrate various behaviors in how the follower maintains the formation position. The developed tool has proven adequate to support future research in crafting real-time controllers capable of determining near-optimal trajectories.  
\end{abstract} 
\tableofcontents

\section{Introduction}
Autonomous formation control is a growing area of interest for future Air Force operations. Autonomous aircraft are increasingly used across the battle space to avoid putting personnel at risk. Thus, the need for autonomous aircraft to handle formation tasks normally requiring a pilot has grown. One model for this interaction consists of a leader aircraft who is accompanied by an autonomous follower aircraft. The ability for the follower aircraft to autonomously rejoin and maintain a designated formation position would free the leader aircraft to focus on other duties, leading to increased mission effectiveness. The leader aircraft often has primary responsibility for interactions outside of the formation, so having an autonomous wingman which can execute without additional oversight from the lead aircraft in an efficient manner would decrease the leader workload. The scenario investigated in this paper focused on the response of a wingman rejoining from an arbitrary starting location to a defined formation position and maintaining that position until a fixed final time. This same behavior could readily be applied to the scenario of an autonomous aircraft attempting to rapidly attain an attack position on an enemy aircraft, along with many other combat and training applications.

\begin{figure}[h]
\centering
\includegraphics[width=3.2in]{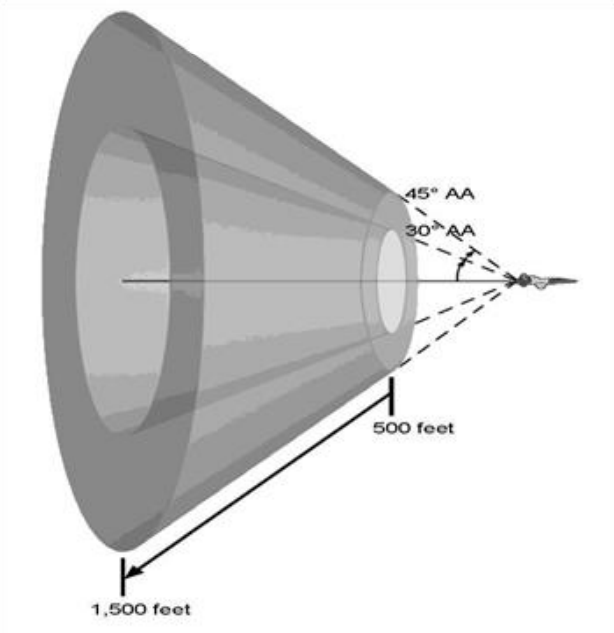}
\caption{\bf{Fighting Wing Position} \protect\cite{a251}}
\label{fw}
\end{figure}

 Extensive research in maintaining a desired formation position, defined as a single location relative to a leader aircraft, has been conducted in articles such as \cite{Buzogany1993}, \cite{Pachter1994}, and \cite{Reyna1994}. The author of \cite{WeintraubRejoin} presented a model predictive controller which added the ability to arrive at a target waypoint with a desired velocity vector orientation. To extend this concept beyond single-point formation positions, further research was conducted which crafted a control architecture allowing a follower aircraft to rejoin to a ring of points defined relative to a leader aircraft \cite{Tran2021}. This ring of formation points was inspired by the traditional formation position "fighting wing" as seen in Figure \ref{fw} referenced from \cite{a251}. Future research aims to dictate a trajectory in real-time for which an autonomous wingman would attain and maintain a formation position defined in reference to a maneuvering leader while minimizing a desired cost function. In order to attain an approximate minimum in real-time may require the use of different control laws depending on the initial states of the leader and follower aircraft. Thus, a state space of boundaries could be defined to dictate when certain control laws would be optimal. Further, future research aims to create a strategy which will optimally select from these predefined control laws and implement them in real-time. Since it is unlikely a fully optimal solution can be attained in real-time for this problem, the control strategy will aim to approximate an optimal trajectory. Therefore, the ability to calculate the degree of the error of this approximation is needed. This paper specifically aimed to create a validation tool that will determine the optimal path for future comparison to the approximated trajectory generated by the real-time control strategy.

One specific scenario was considered. The follower aircraft being tasked to rejoin to a defined formation ring in minimum time, and then maintaining that formation position with minimum deviation over a fixed time duration. During this task, the leader aircraft maneuvered in all three axes, similar to the behavior dictated by typical formation roles. The leader's maneuvers were determined \textit{a priori} for these simulations, and were treated as a given set of parameters for each test case. To simulate the follower aircraft's flight envelope, constraints were applied to the magnitude of both the follower's velocity and acceleration. The leader aircraft maintained its own constraints along the predetermined path. Two cases were analyzed and evaluated to assess the tools effectiveness at determining the optimal path. The analysis of this problem is organized in the following sections: Section \ref{PF} details the problem formulation, Section \ref{SM} outlines the solution methodology, Section \ref{R} includes discussion of the results, and Section \ref{C} includes the conclusions and future recommendations.

\section{Problem Formulation}\label{PF}  

\subsection{Assumptions}
To aid in the initial creation of the validation tool, a number of simplifying assumptions were used. The follower aircraft was modeled with triple integrator dynamics, with control applied as a three-dimensional jerk vector. This model assumed perfect knowledge of the aerodynamic forces and thrust interactions. Full knowledge of the leader's path was provided to the solver for both cases, including three-dimensional position, velocity, acceleration, and jerk vectors defined in the inertial reference frame. No exogenous inputs were considered within the tool. Therefore, deterministic knowledge of all the follower states was assumed, which included three-dimensional position, velocity, and acceleration. 

Both the leader and follower aircraft were modeled as point masses about their center of gravity. For any images created throughout this report, the velocity vector was used as a reference to determine the pitch and the yaw of the aircraft. However, without modeling the full Euler angles of the aircraft the ability to determine the roll angle was limited. To aid in graphics development, roll was visually estimated. 

In formulating the optimal control problem, the existence of an optimal solution was assumed. A scenario without a solution could be easily crafted, such as the simple case where a trailing follower cannot accelerate fast enough to catch a higher performance leader aircraft in the two-minute simulation time. This assumption was satisfied within this study by choosing appropriate initial conditions and using matching aircraft performance envelopes. This assumption could easily be relaxed while still providing optimal solutions in future work. Additionally, the test cases presented leader trajectories such that the opportunity existed for an optimal rejoin and formation position maintenance. 

Lastly, the aircraft performance envelope was modeled with a minimum velocity and acceleration to prevent aerodynamic stall, a maximum velocity to limit maximum dynamic pressure forces, and a maximum acceleration to limit maximum structural forces.

\subsection{Reference Frames}
\begin{figure*}
\centering
\includegraphics[width=5.5in]{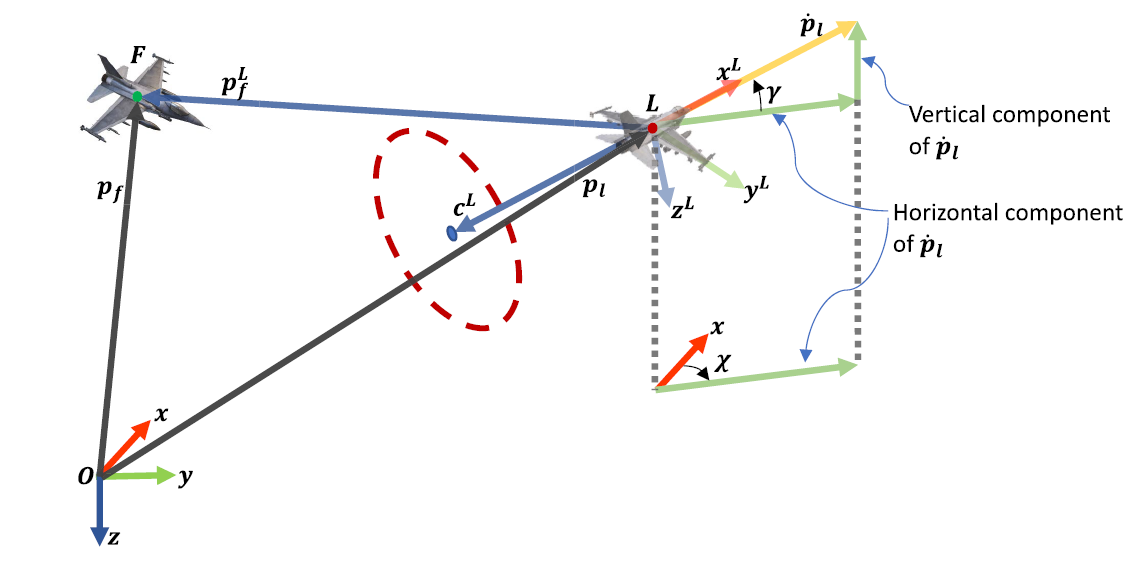}
\caption{\bf{The formation ring is defined in the leader frame as shown by the red circle, in trail of the leader.  The follower strives to rejoin by reaching the ring and then sustaining that position thereafter.} \protect\cite{Tran2021} }
\label{refframe}
\end{figure*}

This problem required two primary reference frames: the inertial reference frame and the leader reference frame. Both are shown in Fig \ref{refframe}, reproduced from \cite{Tran2021}.
The inertial reference frame was defined as a traditional north, east, down right-handed coordinated system with an origin fixed to the starting north and east location of the follower aircraft at zero feet mean sea level (MSL). The body fixed reference frame of interest was fixed to the center of gravity of the leader aircraft. This leader reference frame was defined with a longitudinal x axis for which positive extended out the nose of the aircraft, and a lateral y axis with positive defined as pointing out the right wing. The third axis of the leader frame was defined by the cross product of the x and y axis, using the right-hand rule to determine the sign. The rotation matrix between the two reference frames was defined by the leader's flight path angle, $\gamma$, and the course angle $\chi$. The applicable direction cosine matrix to rotate a vector from the leader frame to the inertial frame was given by the following equation referenced in \cite{Beard2012}.
\begin{equation}
\label{eq:rotation}
    R_L^i = \begin{bmatrix}\cos{ \chi }&-\sin{ \chi}&0\\\sin{\chi}&\cos{ \chi }&0\\0&0&1\end{bmatrix}\begin{bmatrix}\cos{\gamma}&0&\sin{\gamma}\\0&1&0\\-\sin{\gamma}&0&\cos{\gamma}\end{bmatrix}
\end{equation}
Additionally, the rotation matrix in (\ref{eq:rotation}) satisfies the relations $(R_i^L)^{-1}=[R_L^i]^\top=R_L^i$ and $\text{det}(R_i^L)=1$. $R_i^L$ defines the rotation from the inertial frame to the leader frame. A superscript L on any vector denotes a vector defined in the leader reference frame, while all other vectors are defined in the inertial frame. The subscript $l$ and $f$ represent the leader and follower vectors respectively. Lastly, the vector $p_d$ represent the vector \textbf{LF} in the inertial frame as shown in Figure \ref{refframe}. This formulation was chosen to match the research in \cite{Tran2021}. The follower and leader states and control jerk vector are defined in relation to the inertial coordinate system as follows:
\setlength{\tabcolsep}{0 pt}
\begin{table}[h]
    \centering
    \begin{tabularx}{0.49\textwidth}{rXl}
       $p_i=[p_{i_x}, p_{i_y}, p_{i_z}]^\top$;  $LF=p_d=p_f-p_l$ \dotfill & \dotfill& \dotfill Position\\
       $v_i=[v_{i_x}, v_{i_y}, v_{i_z}]^\top$\dotfill & \dotfill& \dotfill Velocity\\
       $a_i=[a_{i_x}, a_{i_y}, a_{i_z}]^\top$\dotfill & \dotfill& \dotfill Acceleration\\
       $u_i=[u_x, u_y, u_z]^\top$\dotfill & \dotfill& \dotfill Control\\ 
    \end{tabularx}
\end{table}
\subsection{Mathematical Formulation}
As stated previously, the goal of this validation tool is to determine the optimal path for a follower aircraft to rejoin to a formation position, defined by the set of points contained on the formation ring in minimum time, and subsequently maintain the formation position with least deviation until a fixed final time. This problem was crafted into a two-phase dynamic optimal control problem based on the GPOPS-II standard problem formulation \cite{Patterson2014}.

The cost function for phase 1 dictates a minimum time to execute the initial rejoin from a predefined starting location to a point contained on the formation ring. The Mayer form of a minimum time cost function was used as the final time of phase 1. The follower's dynamics were modeled by three-dimensional triple integrator dynamics, with the three-element control applied directly to the follower's jerk vector. The scenario begins at time zero with a known initial state of both the leader and follower. A path constraint was enforced to prevent the follower from entering the leader's jet wash, a safety hazard to the follower aircraft. The jet wash was modeled as an infinite cylinder centered about the $x^L$-axis centered at $(y^L,z^L)=(0,0)$ with a radius of $R_{jw}$. Note that $t_f^{(i)}$ indicates the final time of phase $i$. Thus, phase 1 included the terminal condition requiring the follower to attain a formation position at $t_f^{(1)}$ which was enforced by the functions:
\begin{equation}
\begin{tabular}{l l}
    \small $f_1(p_d) = (p_{d_x}-[R_L^i c^L]_x)^2$\\
    \small $f_2(p_d) = (p_{d_y}-[R_L^i c^L]_y)^2+(p_{d_z}-[R_L^i c^L]_z)^2-R^2_{ring}$
\end{tabular}
\label{eq:ringEq}
\end{equation}
%

The equations in (\ref{eq:ringEq}) provided a way to characterize the distance from the follower aircraft to the formation ring. $f_1(p_d(t_f^{(1)}))$ ensured the follower aircraft attained the same $x^L$-coordinate in the leader reference frame as the tip of the vector $c^L$, which was used to define the center point of the formation ring.  $f_2(p_d(t_f^{(1)}))$ was used to ensure that the $y^L$ and $z^L$ coordinates of the follower in relation to the leader were at a distance from the ring center equal to the radius of the ring. Simply put, when $f_1=f_2=0$ the follower was established in a valid formation position defined by the ring of points of radius $R_{ring}$ with center defined by the vector $c^L$. In order to permit small numerical errors, a tolerance of $\pm 10$ meters was allowed as the maximum value of each function at the terminal condition. Lastly, constraints were enforced on the magnitude of the follower's velocity and acceleration vectors for the duration of the phase. The magnitude of the follower's velocity vector was constrained by a minimum and maximum velocity, as dictated by the aircraft's stall speed and maximum speed. The magnitude of the acceleration vector was constrained by a minimum acceleration to avoid zero g flight, and a maximum acceleration to prevent overloading the aircraft structure. Both of these constraints were enforced as path constraints throughout the duration of the phase.   

Minimizing the cost function of phase 2 required the follower aircraft to minimize the deviation from the formation ring for the duration of the phase. This deviation was characterized by squaring and summing $f_1$ and $f_2$ as defined during phase 1, and integrating the resulting values over the duration of the phase. Phase 2 used the same follower dynamics along with the same path, state, and control constraints. The second phase also contained the requirement for the position and velocity to match across phases with a tolerance of $\pm 0.1$ m, ensuring a continuous trajectory throughout the entire scenario. Additionally, the condition for the free final time of phase 1 to match the initial time of phase 2 was enforced. The final time of phase 2 was fixed, and was chosen to ensure that the follower could attain the formation ring during phase 1 before the fixed final time of phase 2. 

Since the optimization software required one cost function for the overall problem, the cost functions from both phase 1 and phase 2 were summed to produce a cumulative cost function. In order to account for a difference in the magnitudes between the cost functions of each phase, a weighting parameter $\beta$ was included in the summation.
\setlength{\tabcolsep}{0 pt}
 \renewcommand{\arraystretch}{1.5}
\begin{table*}[!ht]
    \centering
    \begin{tabularx}{0.85\textwidth}{rclXr}
       Overall: \ \ \ &$\min\limits_{u(t)}$\ \ \ &  $J=\beta J_1+(1-\beta)J_2$ \dotfill&\dotfill&\dotfill Cost Function\\
       \\
       Phase 1: \ \ \ &$\min\limits_{u(t)}$\ \ \ & $J_1=t_f^{(1)}$\dotfill&\dotfill&\dotfill Cost Function\\
       &S.T. \ \ \ &$\dot{p_f}=v_f$\dotfill&\dotfill&\dotfill Dynamics\\
       & & $\dot{v_f}=a_f$\dotfill&\dotfill&\dotfill Dynamics\\
       & & $\dot{a_f}=u$\dotfill&\dotfill&\dotfill Dynamics\\
       & & $t_0=0\text{, }p_f(t_0)=p_{f_0}\text{, }v_f(t_0)=v_{f_0}\text{, }a_f(t_0)=a_{f_0}$\dotfill&\dotfill&\dotfill Boundary Condition\\
       & & $-\epsilon_{form}\leq f_1(p_d(t_f^{(1)}))=(p_{d_x}(t_f^{(1)})-[R_L^i(\gamma, \chi)c^L]_x)^2\leq \epsilon_{form}$\dotfill&\dotfill&\dotfill Terminal Condition\\
       & & $-\epsilon_{form}\leq f_2(p_d(t_f^{(1)}))=(p_{d_y}(t_f^{(1)})-[R_L^i(\gamma, \chi)c^L]_y)^2$& & \\
       & & $\ \ \ \ \ \ \ \ \ \ \ \ \ \ \ \ \ \ \ \ +(p_{d_z}(t_f^{(1)})-[R_L^i(\gamma, \chi)c^L]_z)^2-R_{ring}^2\leq \epsilon_{form}$\dotfill&\dotfill&\dotfill Terminal Condition\\
       & & $-([R_i^L(\gamma, \chi)p_d]_y)^2-([R_i^L(\gamma, \chi)p_d]_z)^2+R_{jw}^2 \leq 0$\dotfill&\dotfill&\dotfill Jet Wash\\
       & & $v_{f_{min}}\leq |v_f|\leq v_{f_{max}}\text{, }a_{f_{min}}\leq |a_f|\leq a_{f_{max}}$\dotfill&\dotfill&\dotfill Constraints\\
       \\
       Phase 2: \ \ \ &$\min\limits_{u(t)}$\ \ \ & $J_2=\int\limits_{t_0^{(2)}}^{t_{f}^{(2)}}\{f_1(p_d(t))^2+f_2(p_d(t))^2\}dt$\dotfill&\dotfill&\dotfill Cost Function\\
       &S.T. \ \ \ &$\dot{p_f}=v_f$\dotfill&\dotfill&\dotfill Dynamics\\
       & & $\dot{v_f}=a_f$\dotfill&\dotfill&\dotfill Dynamics\\
       & & $\dot{a_f}=u$\dotfill&\dotfill&\dotfill Dynamics\\
       & & $t_f^{(2)}=120$\dotfill&\dotfill&\dotfill Boundary Condition\\
       & & $t_o^{(2)}=t_f^{(1)}\text{, }p_f(t_0^{(2)})=p_f(t_f^{(1)})$\dotfill&\dotfill&\dotfill Phase Continuity\\
       & & $v_f(t_0^{(2)})=v_f(t_f^{(1)})\text{, }a_f(t_0^{(2)})=a_f(t_f^{(1)})$\dotfill&\dotfill&\dotfill Phase Continuity\\
       & & $-([R_i^L(\gamma, \chi)p_d]_y)^2-([R_i^L(\gamma, \chi)p_d]_z)^2+R_{jw}^2 \leq 0$\dotfill&\dotfill&\dotfill Jet Wash\\
       & & $v_{f_{min}}\leq |v_f|\leq v_{f_{max}}\text{, }a_{f_{min}}\leq |a_f|\leq a_{f_{max}}$\dotfill&\dotfill&\dotfill Constraints\\
    \end{tabularx}
\end{table*}
\setlength{\tabcolsep}{6 pt}
 \renewcommand{\arraystretch}{1}
\setlength{\tabcolsep}{6 pt}
 \renewcommand{\arraystretch}{1.5}
 \begin{table*}[tp]
 \centering
\caption{Initial Conditions}
\label{ictable}
\begin{tabular}{|l|l|l|l|}
\hline
$p_{f_0}=\begin{bmatrix}0 \text{ ft}\\0 \text{ ft}\\20,000 \text{ ft}\end{bmatrix}$ & $v_{f_0}=\begin{bmatrix}450 \text{ kts}\\0 \text{ kts}\\0 \text{ kts}\end{bmatrix}$ & $a_{f_0}=\begin{bmatrix}0.07 \text{ g}\\0 \text{ g}\\1 \text{ g}\end{bmatrix}$ & $c^L=\begin{bmatrix}-700 ft\\0\\0\end{bmatrix}$ \\ \hline
$R_{ring}$ = 700 ft & $R_{jw}$=5 ft & $\beta = 0.9$ & $\epsilon _{form} = 10$ m \\ \hline
$V_{min}$ = 200 kts & $V_{max}$ = 700 kts & $A_{min}$ = 0.25 g & $A_{max}$ = 7 g\\ \hline
\end{tabular}
\end{table*}

A summary of the problem's mathematical formulation is presented on the following page.

\subsection{Discrete Hamiltonian Derivation}
A derivation of the discrete Hamiltonian is provided and was used to ensure the validity of the optimal solution. The continuous Hamiltonian was defined as $\mathcal{H}=\mathcal{L}+\Bar{\lambda}^\text{T}\mathbf{F}$. The Lagrangian of the cost function was defined as $\mathcal{L}$$=(1-\beta )\{f_1^2+f_2^2\}$. $\Bar{\lambda}$ contained a column vector of the 9 costates. Finally, $\mathbf{F}$ was a column vector of the right-hand side of the 9 dynamics equations. Since the tool produces a solution at discrete points in time, the discrete form of the Hamiltonian equation was used, defined as $\mathcal{H}$$_d=\mathcal{L}$$_k+\Bar{\lambda}_k^\text{T}\mathbf{F}_k$. Here the subscript $k$ denotes the value of each function at the point $t_k$. For a problem in which the Hamiltonian is not an explicit function of time, the value of the Hamiltonian will be constant on an extremal path \cite{kirk}. This criterion is applicable to both phase 1 and phase 2. 

\subsection{Test Cases}
 A summary of the initial conditions and parameters is presented in Table \ref{ictable}. Both aircraft were modeled as a typical 7g capable fighter aircraft. The formation ring center was set at 700 ft directly behind the leader aircraft. The radius of the formation ring was set to 700 ft. This placed the formation ring points approximately 1,000 ft from the leader aircraft on a 45-degree aspect, directly correlating to the fighting wing position discussed previously. The jet wash cylinder was set 5 feet wide. The follower aircraft was initialized in 1 NM trail of the leader aircraft, with a 1,000 ft lateral offset to the leader's right wing. The simulation was run for 2 minutes. The velocity was limited to a minimum of 200 knots to prevent stall, and a maximum of 700 knots. The maximum acceleration was limited to 7 g's, with a minimum of 0.25 g's.
\begin{figure*}[bp]
\centering
\includegraphics[width=5.98in]{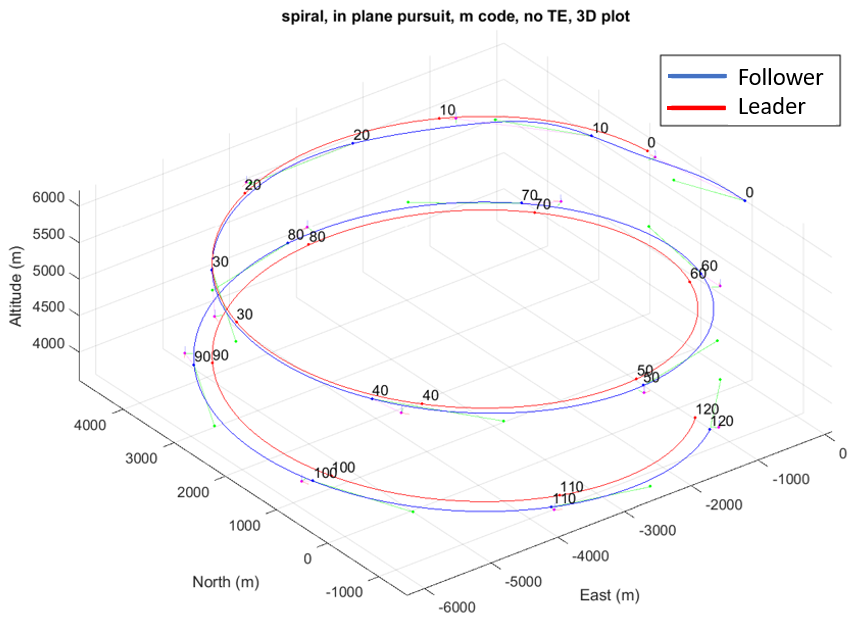}\\
\caption{\textbf{Initial Guess Provided by Algorithm}}
\label{icfig}
\end{figure*}

 In order to assess the capabilities of the validation tool, two test cases were evaluated. Both aircraft began on a north heading in level flight at 450 knots and 20,000 ft MSL. The first leader trajectory consisted of a constant left bank, descending spiral. The leader descended with constant airspeed at approximately 3.25 degrees of pitch losing 6,500 ft over the 2-minute simulation. The second case used the same initial conditions followed by the leader performing a number of continuous loops. The loops were approximately 12,000 ft tall starting and ending at the initial altitude. The loops were executed on a constant heading. Again, a 2-minute simulation time was used. Due to the minimum time formulation of the problem, the optimal solution was expected to take a direct path to the formation ring, with allowances for maintaining the ring position throughout phase 2.

\section{Solution Methodology}\label{SM}
\subsection{Direct Collocation Solver}
This problem was solved with the direct orthogonal collocation solver GPOPS-II discussed in \cite{Patterson2014}. GPOPS-II uses a polynomial approximation of the state and control vectors. The solver transcribed the continuous optimal control problem by interpolating these approximations at LGR based collocation points as discussed in \cite{Garg2010}. Gaussian quadrature was used to approximate the cost function, and a differential matrix based on the Lagrange polynomial basis allowed the dynamics to be enforced exactly at each collocation point. Additionally, the boundary and path constraints were enforced at each collocation point. GPOPS-II then formatted the problem for hand-off to a nonlinear program (NLP) solver. The NLP solver used both the gradient and Hessian of the transcribed problem. These functions were computed using the MATLAB program ADiGator discussed in \cite{weinstein2017}. 

Scaling was used to aid in solver convergence, which simply consisted of normalizing to 100s of meters of distance. The NLP solver requires bounds to constrain the feasible region for all variables. The position state was bounded to the max distance the follower aircraft could travel straight along the coordinate direction at maximum velocity over the two-minute simulation time. The velocity and acceleration components were bounded by $\pm V_{max}$ and $\pm A_{max}$ respectively. An acceleration rate limit of $\pm 5$ g/s was used to bound the control components. Unless specified elsewhere, default settings were used as provided with GPOPS-II.   

Once the first solution was returned from the NLP solver, the error in the enforcement of the dynamics at the midpoints between collocation points was estimated by GPOPS-II. The maximum magnitude of this error was used as the termination criteria for the algorithm, with a tolerance of $10^{-3}$. If the error exceeded the tolerance, GPOPS-II modified the mesh of each phase, and the number of collocation points included on each mesh interval. The hp-LiuRao adaptive mesh refinement method provided with GPOPS-II was used. This process was then repeated until the mesh tolerance was satisfied. GPOPS-II reported the final solution with values of the states, controls, and costate estimates at each collocation point and one interpolated point at the end of each of the two phases.  

\subsection{NLP Solver}
The NLP Solver IPOPT was used as provided with GPOPS-II. IPOPT used an interior point boundary value method which was described in \cite{wachter}. A tolerance of $10^{-5}$ was used for the determination of convergence within IPOPT. The linear solver MA57 was used by IPOPT to solve the system of equations on each NLP iteration. 

\subsection{Initial Guess Generation}

GPOPS-II required an initial guess for the states and controls over the duration of the problem. In order to align with future work, the initial guess of the follower's states and controls were generated from a control strategy algorithm. This algorithm was currently under development with the intent of creating the ability to approximate the optimal trajectory in real-time. Details of this algorithm will be documented in a future publication. The algorithm used a control strategy which selected from simple maneuver primitives based on various decision criteria informed by the relation between the leader and follower aircraft at a given instance in time. While this algorithm is still under development, it has a basic capability such that the follower aircraft can track to a single point defined in reference to the leader aircraft.

Using this basic capability, the algorithm was run such that the follower tracked to and maintained a point that would coincide with a single point located on the formation ring. This algorithm did not allow for the tracking of a set of points at the time of this writing, as the validation tool does. For both cases the initial guess was generated by having the follower aircraft attempt to track to the point that is longitudinally in line with the projection of the right-wing tip of the leader. While this algorithm was still under development, it provided a sufficient initial guess such that the optimal control solver was able to converge to an optimal trajectory in all cases. 

In order to assess the impact various initial guesses may have had on the convergence of the optimal control solver, the spiral case was evaluated with 4 different algorithm results. The algorithm was set to track to 4 various points along the formation ring. The four points were located along the projection of the left wing, along the projection of the right wing, and at the top and bottom of the ring when referencing the $z^L$ axis. The algorithm results were then evaluated on their time to rejoin to the formation ring and the magnitude of the final cost function for comparison. An example of an initial guess for the spiral trajectory can be seen in figure \ref{icfig}.
\begin{figure*}[bp]
\centering
\includegraphics[width=5.98in]{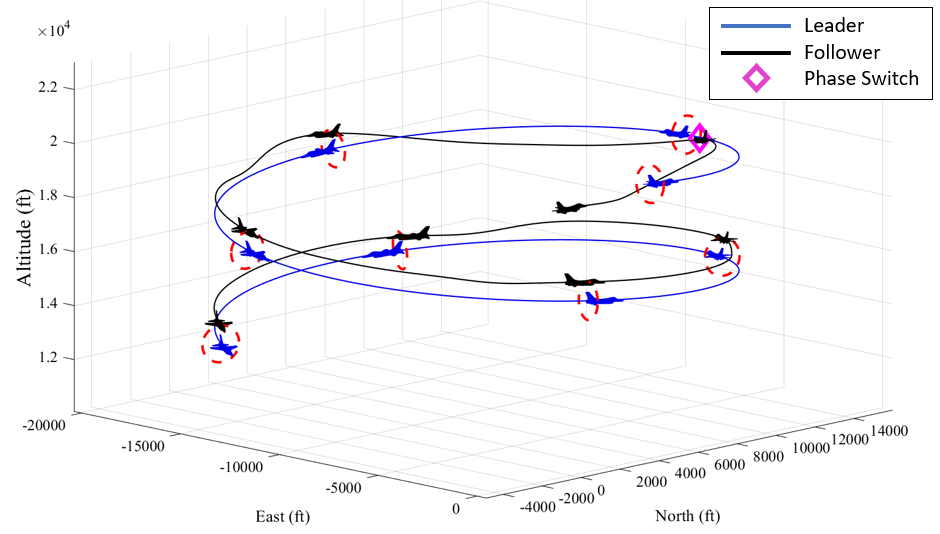}
\caption{\bf{Position Response for Case 1}}
\label{posc1}
\end{figure*}

\section{Results}\label{R}
\subsection{Spiral Trajectory}
\begin{figure}[!ht]
\centering
\includegraphics[width=3.2in]{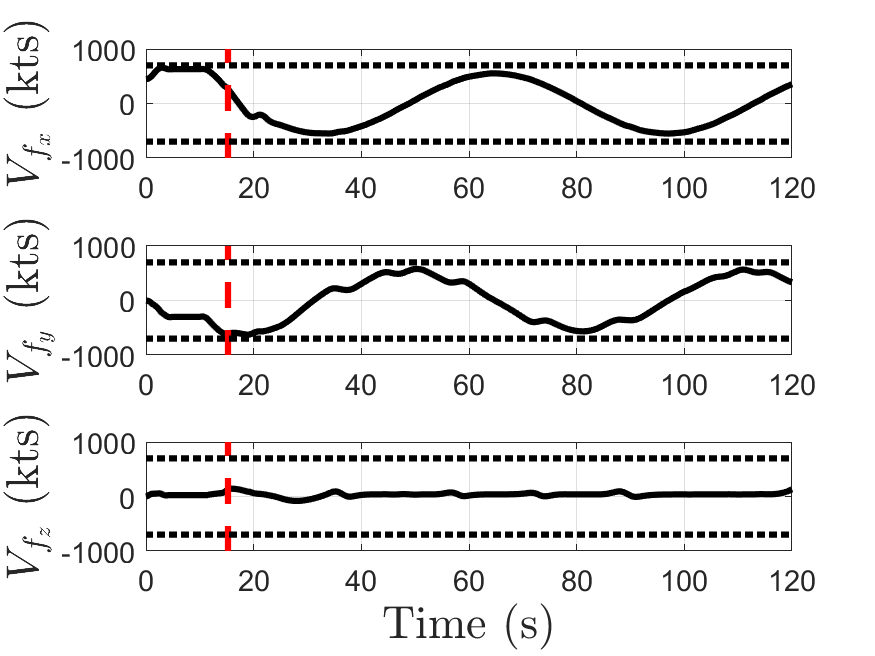}
\caption{\bf{Velocity Response for Case 1}}
\label{velc1}
\end{figure}

\begin{figure}[!ht]
\centering
\includegraphics[width=3.2in]{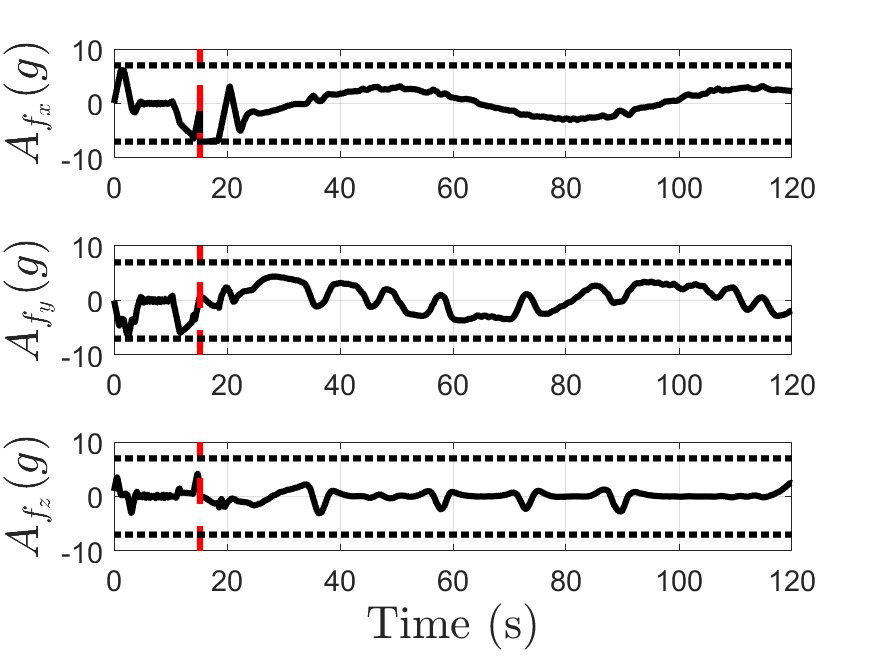}
\caption{\bf{Acceleration Response for Case 1}}
\label{accelc1}
\end{figure}

\begin{figure}[!ht]
\centering
\includegraphics[width=3.2in]{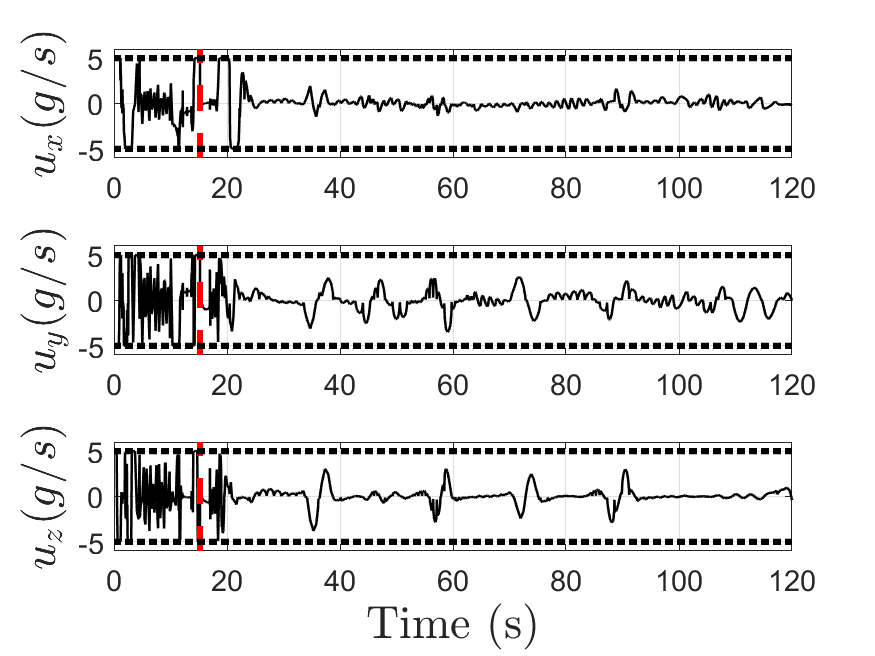}
\caption{\bf{Optimal Control for Case 1}}
\label{controlc1}
\end{figure}

For the descending spiral trajectory, GPOPS-II reported an optimal solution for which the states and controls are depicted in Figures \ref{posc1}, \ref{velc1}, \ref{accelc1}, and \ref{controlc1}. The position response showed the follower aircraft in black, the leader aircraft in blue, and a maroon diamond represented the point in the trajectory when the phase transition occurs. As expected, the follower took a direct path to the formation ring during phase 1, while maintaining the formation ring nearly perfectly throughout the duration of phase 2. It should be noted that some motion along the formation ring was observed throughout phase 2, but as this incurs no penalty within the cost function this does not preclude this path from being optimal. The velocity and acceleration response showed the aircraft attempting to achieve a maximum velocity in order to minimize the rejoin time. The limits on the jerk control vector could also be seen by the linear changes in the acceleration components. The optimal control is very noisy during phase 1. However, this was the main reason to use triple integrator dynamics as this noise simply represents a g onset rate which wouldn't be prohibitively difficult to implement on a real aircraft. 

\begin{figure}[!ht]
\centering
\includegraphics[width=3.2in]{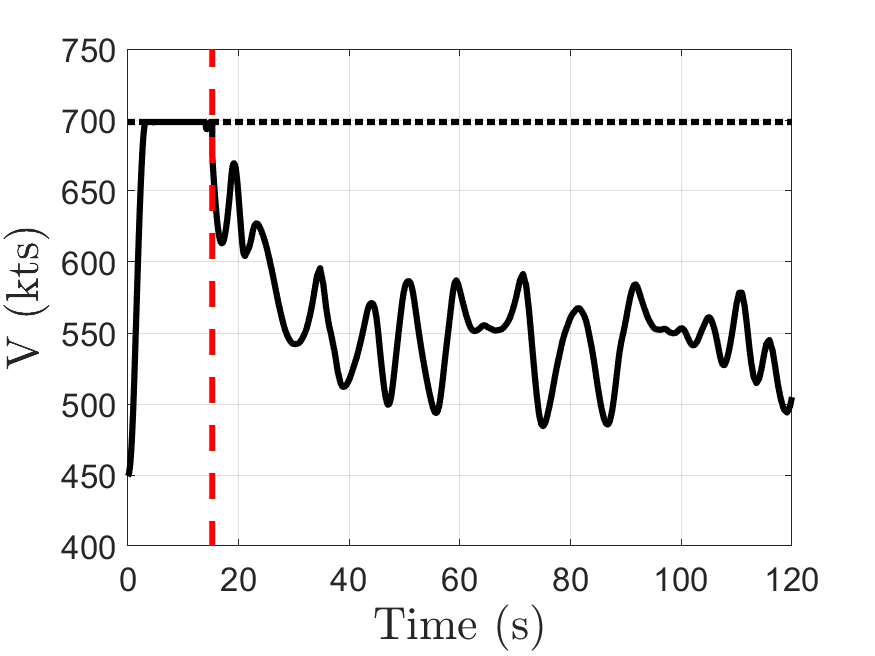}
\caption{\bf{Magnitude of Velocity for Case 1}}
\label{totVc1}
\end{figure}

\begin{figure}[!ht]
\centering
\includegraphics[width=3.2in]{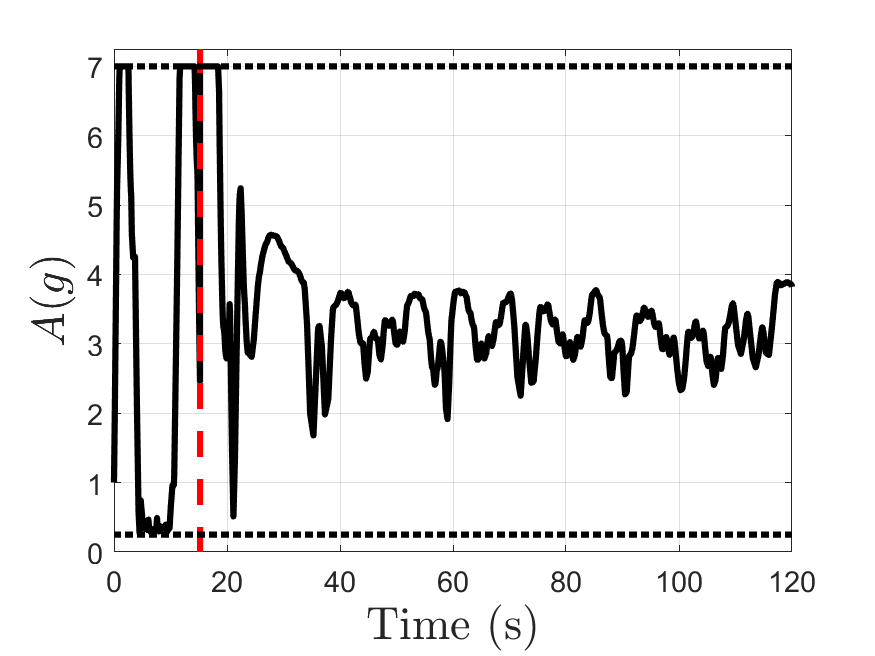}
\caption{\bf{Magnitude of Acceleration for Case 1}}
\label{totAc1}
\end{figure}

The path constraints imposed on the magnitudes of velocity and acceleration can be seen in Figures \ref{totVc1} and \ref{totAc1}. These plots showed the follower aircraft achieving its maximum velocity at the maximum acceleration rate throughout phase 1, reducing the rejoin time to the greatest extent possible. The rejoin time reduction was further enabled by the ability of the follower to maneuvering along the formation ring and dissipate excess energy during phase 2 with no additional cost penalty. This reduction in rejoin time would not have been possible if the formation position was represented by a single point as the follower would have to arrive at that point with its velocity already matched to the leader's velocity in order to minimize deviation during phase 2. \pagebreak

\ \\

\begin{figure}[!ht]
\centering
\includegraphics[width=3.2in]{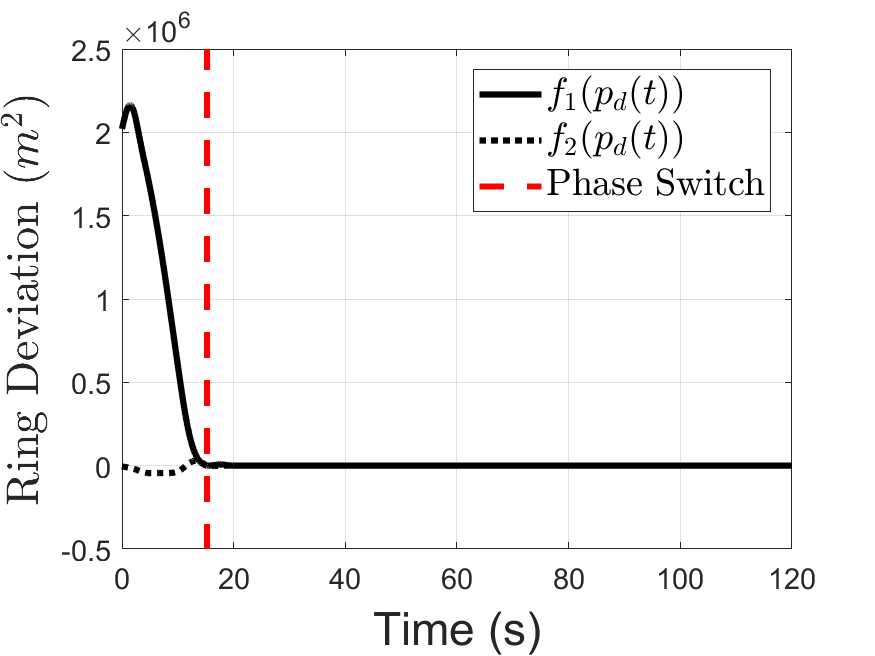}
\caption{\bf{Deviation from Formation Ring in Case 1}}
\label{deviationc1}
\end{figure}

The final plot seen in figure \ref{deviationc1} for case 1 shows the progression of the ring deviation functions throughout the simulation. As expected, the follower rapidly reduced the distance to zero during phase 1, and maintained a zero deviation throughout phase 2 even as the leader continued to maneuver.

\begin{figure*}[bp]
\centering
\includegraphics[width=6.48in]{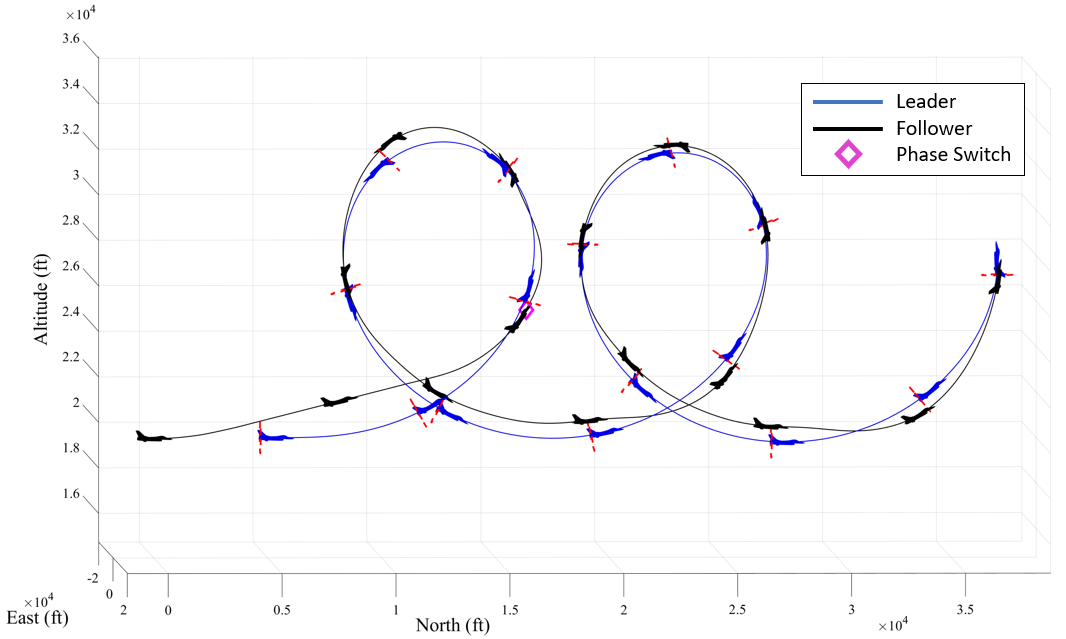}
\caption{\bf{Position Response for Case 2}}
\label{posc2}
\end{figure*}

\subsection{Loop Trajectory}

For the leader trajectory consisting of continuous vertical loops, the optimal path can be seen in Fig \ref{posc2}. This leader trajectory was chosen for study due to its significant differences from the descending spiral. The leader trajectory contained significant vertical deviations along with large changes in the leader's velocity. However, similar behavior of the optimal trajectory was seen across the two cases. During phase 1, the follower again accelerated at the maximum rate to maximum velocity in a direct path to the formation ring. The follower was then able to maintain the ring with little deviation for the duration of phase 2. 

\begin{figure}[!ht]
\centering
\includegraphics[width=3.2in]{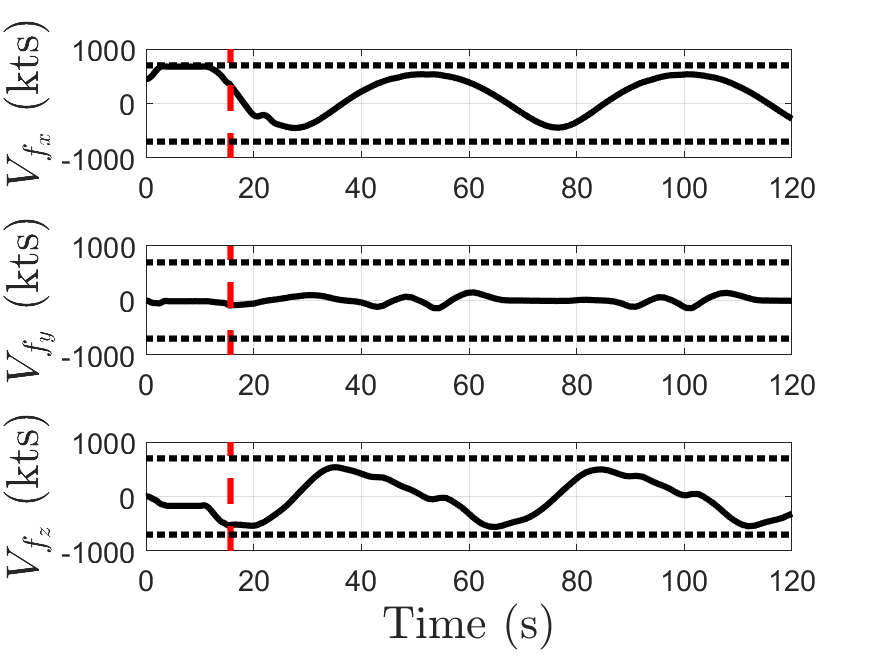}
\caption{\bf{Velocity Response for Case 2}}
\label{velc2}
\end{figure}

\begin{figure}[!ht]
\centering
\includegraphics[width=3.2in]{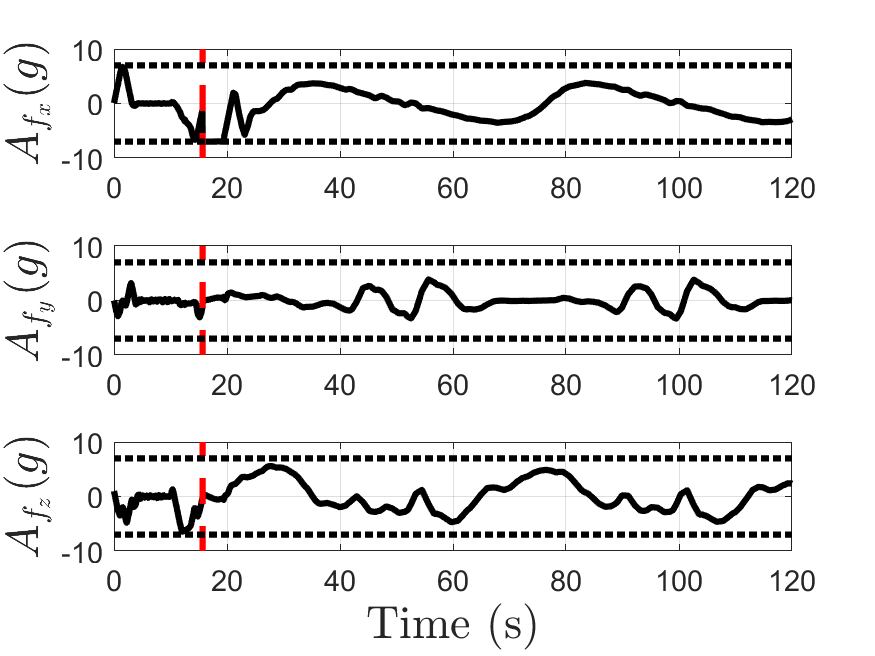}
\caption{\bf{Acceleration Response for Case 2}}
\label{accelc2}
\end{figure}

\begin{figure}[!ht]
\centering
\includegraphics[width=3.2in]{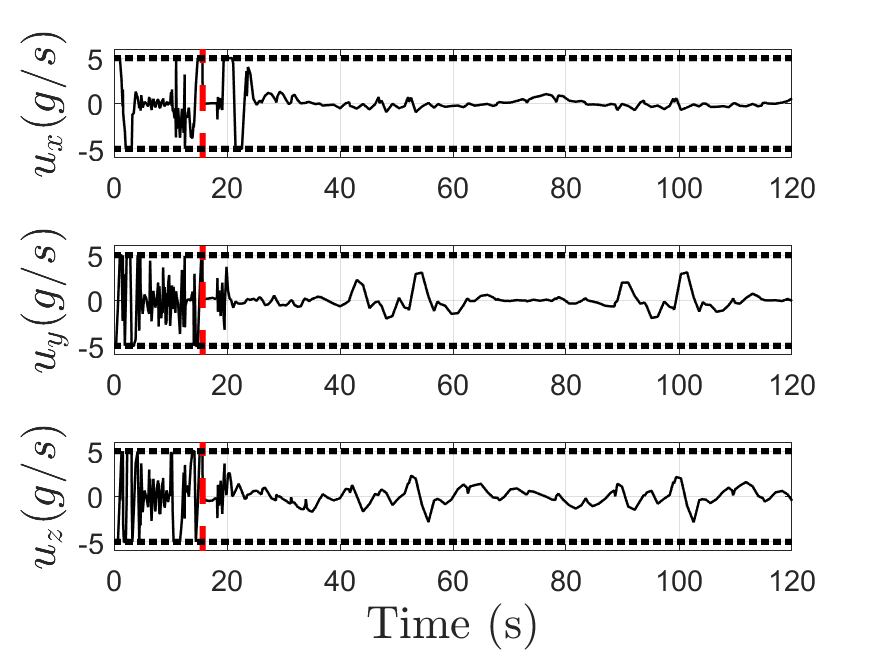}
\caption{\bf{Optimal Control for Case 2}}
\label{controlc2}
\end{figure}

\begin{figure}[!ht]
\centering
\includegraphics[width=3.2in]{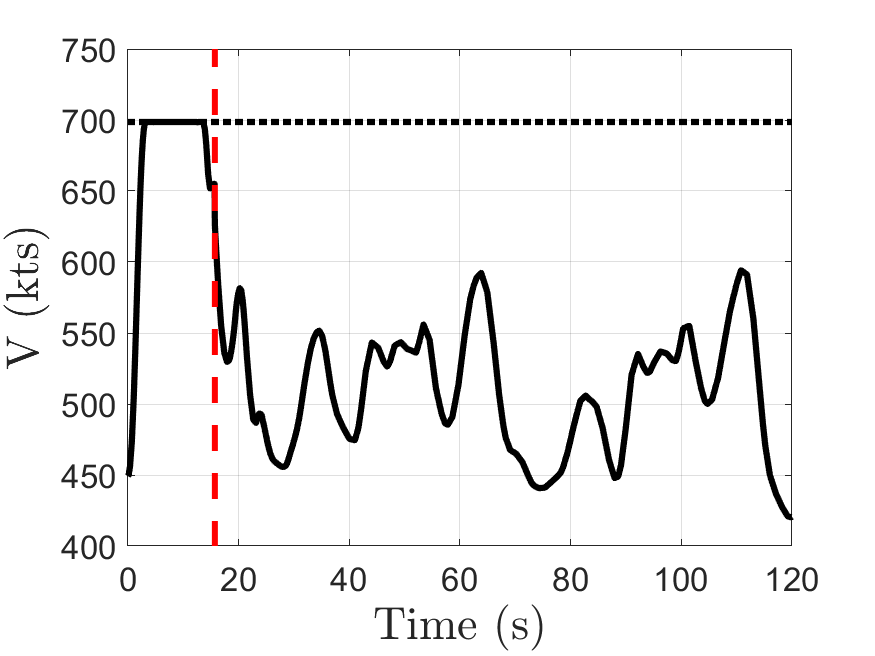}
\caption{\bf{Magnitude of Velocity for Case 2}}
\label{totVc2}
\end{figure}

\begin{figure}[!ht]
\centering
\includegraphics[width=3.2in]{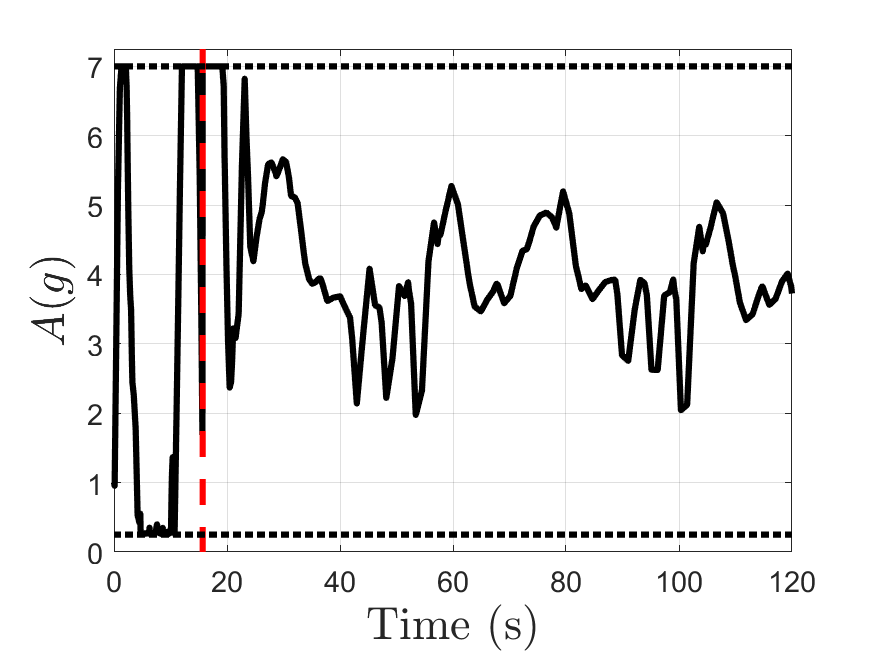}
\caption{\bf{Magnitude of Acceleration for Case 2}}
\label{totAc2}
\end{figure}

\begin{figure}[!ht]
\centering
\includegraphics[width=3.2in]{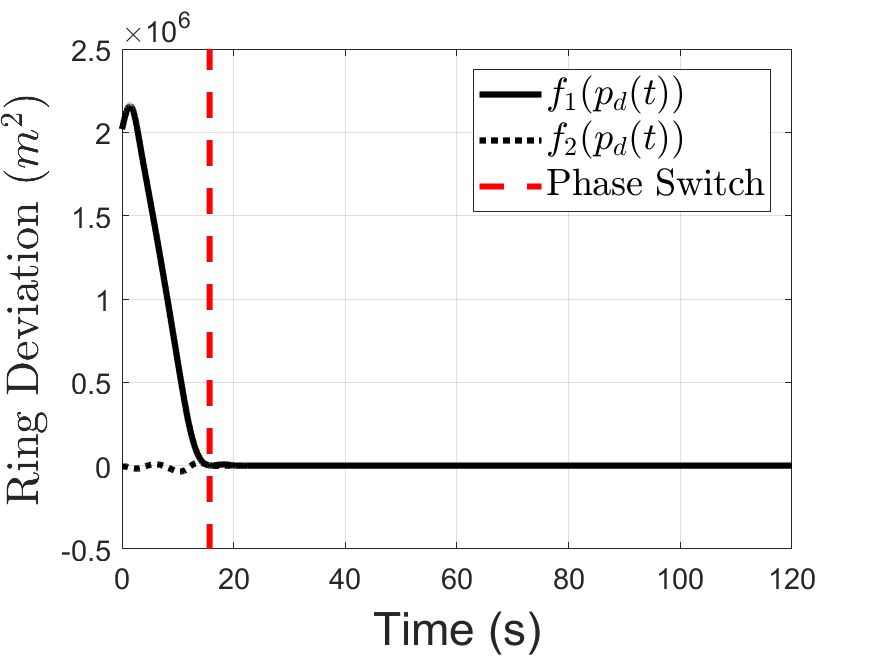}
\caption{\bf{Deviation from Formation Ring in Case 2}}
\label{deviationc2}
\end{figure}

Figures \ref{velc2}-\ref{totAc2} show the other state responses for the follower's trajectory. The maximum velocity attained during phase 1 can be clearly seen, and the large variations in the vertical direction match with the loop maneuvers of the leader. Acceleration rate limiting was also observed in the same fashion as the spiral trajectory. The follower still maneuvered along the formation ring during phase 2 while dissipating excess energy from the rejoin, but with little deviation, as can be seen in Figure \ref{deviationc2}. 

\renewcommand{\arraystretch}{1.7}
\begin{table}[tp]
\centering
\caption{Initial Guess Cost Comparison}
\label{iccomparetable}
\begin{tabular}{|c|c|c|c|c|}
\hline
 & Left Wing & Right Wing & Top & Bottom \\ \hline
$t_f^{(1)}$ & 15.2319 s & 15.2320 s & 15.2422 s & 15.2362 s\\ \hline
$J$ & 16.1472 & 16.1471 & 16.1221 & 16.1221\\ \hline
\end{tabular}
\end{table}

\subsection{Initial Guess Comparison}
\begin{figure*}[tp]
\centering
\includegraphics[width=5.5in]{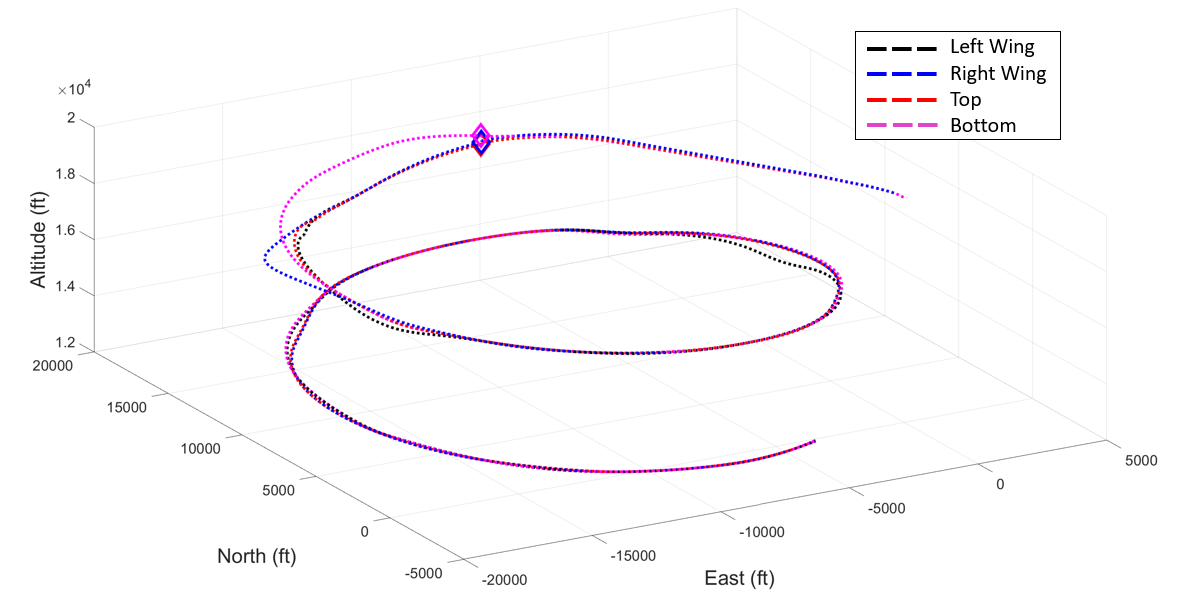}
\caption{\bf{Initial Condition Trajectory Comparisons}}
\label{ICcompare}
\end{figure*}

The second phase of this problem allows for motion along the formation ring without penalty. Thus, more than one path could be considered equally optimal during the second phase. A study of this behavior was conducted by varying the initial guess given to GPOPS-II. Figure \ref{ICcompare} shows the resulting trajectories that were reported as optimal by the software for the 4 tested initial conditions. As predicted, the trajectories were nearly identical during phase one as they took a direct path to minimize the time of the rejoin. During the later portion of phase 2 the trajectories all settled to stable locations along the ring with only small deviations observed between the paths. However, when the follower was dissipating excess energy carried from the first phase, the paths diverged. 

\begin{figure}
\centering
\includegraphics[width=3.2in]{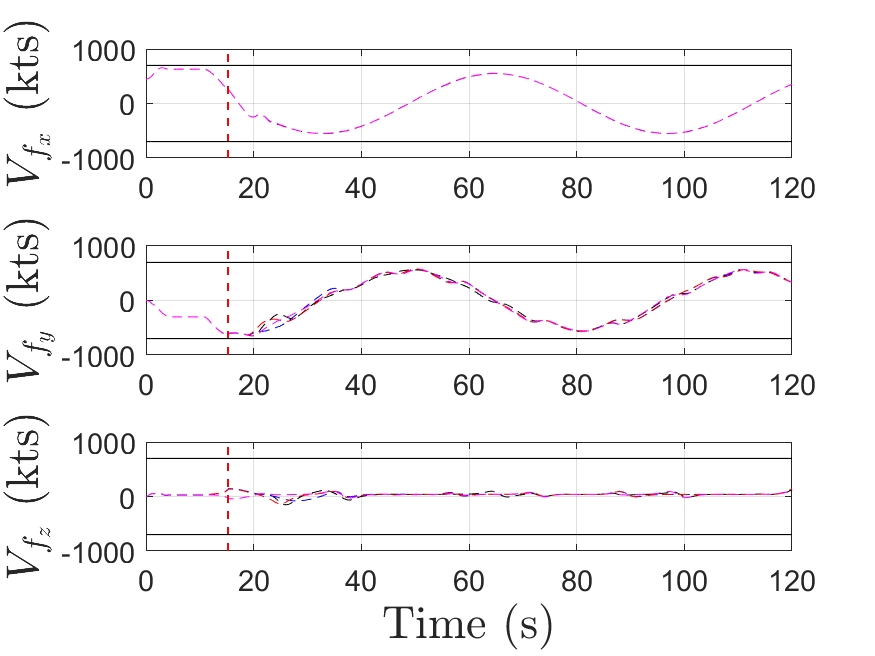}
\caption{\bf{Initial Condition Velocity Comparisons}}
\label{ICcompare_vel}
\end{figure}

Table \ref{iccomparetable} shows the final times of phase 1 and overall cost function values for each initial condition. It is notable that the cost functions are nearly identical even with the path deviations that occurred at the beginning of phase 2. While the various initial conditions caused the solver to converge to somewhat different phase 2 paths, the solver still converged to paths that were equally optimal in regards to the cost function that was studied. This result brings interest to future studies which could include further penalties in the cost function such as limiting control usage or energy expenditure. Figure \ref{ICcompare_vel} shows the similarity in the velocity profiles between each initial condition. This figure uses the same color scheme as Figure \ref{ICcompare}. Similar results were also seen for the other states and the control vector. 
 \ \\

\pagebreak[4]

\hfill\\
\hfill\\
\hfill\\

In all cases the tool produced feasible trajectories with no violation of the path constraints. This was observed for the jet wash constraint by noting that the value of the constraint never took on a positive value. However, the maximum value of the jet wash constraint, $-8.54 \times 10^{-5}$, occurred during the spiral trajectory with the goal position at the top of the ring  So, the optimal trajectory did closely approach the jet wash cylinder, but never violated the path constraint. 

\begin{figure}[!ht]
\centering
\includegraphics[width=3.2in]{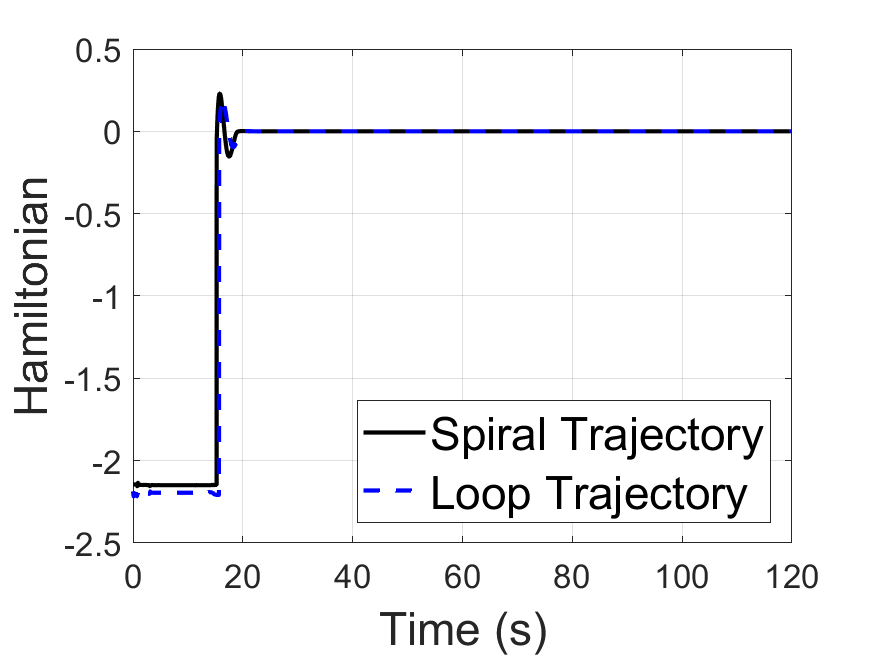}
\caption{\bf{Discrete Hamiltonian}}
\label{hamc1}
\end{figure}

\subsection{Analysis of Optimally}
The discrete Hamiltonian is presented for both leader trajectories in Figure \ref{hamc1}. As previously discussed, the Hamiltonian should be constant during phase 1 and 2. Both figures show that using the costate estimates provided by GPOPS-II did yield Hamiltonians which were consistent with these conditions. There was some variation near the phase switch likely caused by numerical error. However, the magnitude of this error was not observed to have significantly affected the validity of the solutions. 

These results demonstrate the initial ability of the tool to produce optimal trajectories for a variety of complex trajectories.

\section{Conclusions}\label{C}
This paper demonstrated the initial capability of a validation tool which can find the optimal path for a follower aircraft to rejoin to a set of dynamic formation points in minimum time and maintain that position for a given leader trajectory. Of note was the ability of the follower to use the increased flexibility of the formation ring to decrease the rejoin time during the first phase with no penalty during the second phase of the trajectory. It allows for a follower to maintain maximum flexibility while maintaining formation position. This in turn allows the leader aircraft freedom to perform aggressive maneuvering if required. Additionally, the solver showed consistent behavior for two drastically different leader trajectories. The optimal trajectory worked to quickly rejoin to the formation ring, and maintained the ring with nearly zero deviation during phase 2, regardless of the leader maneuvering laterally or vertically. Additionally, the validity of the solutions was shown through the consistent behavior of the discrete Hamiltonian. This tool should prove useful in future research for finding the optimal path to evaluate the accuracy of a control strategy. 

Creating this validation tool uncovered a number of lessons learned for future use. The tool produced a noisy control response, which necessitated the use of triple integrator dynamics to yield realistic acceleration and velocity vectors. While choosing different initial conditions did produce somewhat different optimal trajectories, they did not differ significantly in their rejoin time or total amount of deviation from the formation ring. This result demonstrated that there was a family of optimal solutions which were all approximately equal in value for the given cost function. Future research could consider adding additional cost parameters to penalize things such as excess control effort or energy expenditure.   

A number of areas are planned for future research. The first area was planned to test additional leader trajectories to continue building trust in the tool's ability to generate optimal trajectories. A more complicated dynamics model will be incorporated for the follower's dynamics, incorporating aircraft aerodynamics. Different formation shapes can be defined by modifying the cost functions accordingly. However, the next major area of research will consist of determining control laws and strategies which can approximate these optimal trajectories, but can be computed at a speed which will allow for real-time implementation.

\acknowledgments
The authors thank the Air Force Research Laboratories for funding this project.

\bibliographystyle{IEEEtran}
\bibliography{library}

\thebiography
\begin{biographywithpic}
{Carl Gotwald}{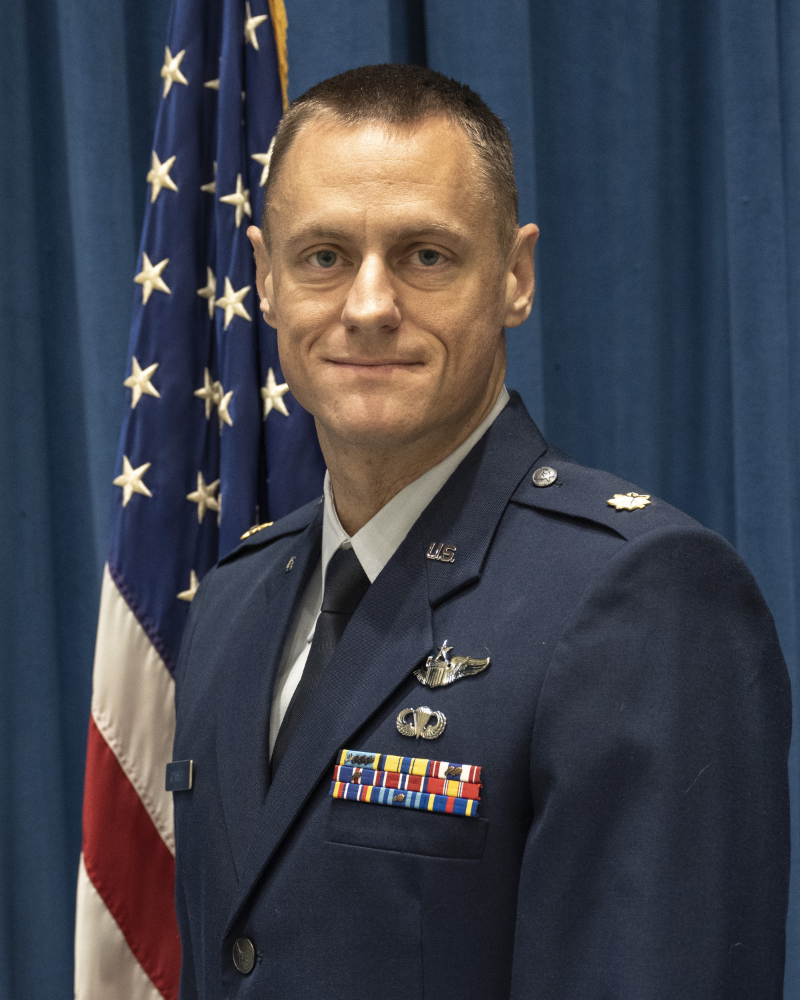}
received his B.S. from the United States Air Force Academy. He holds M.S. degrees from both Mississippi State University and Air University, and is currently pursuing a Ph.D. at the Air Force Institute of Technology. He works as a test pilot with extensive experience testing remotely piloted aircraft and instructing at the United States Air Force Test Pilot School. His current research interests lie in optimal control and real-time simulation. 
\end{biographywithpic} 

\begin{biographywithpic}
{Michael Zollars}{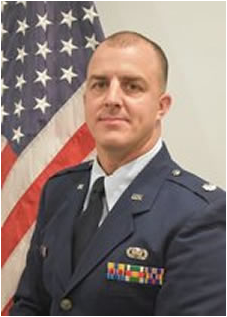}
received a B.S. in Mechanical Engineering from the Pennsylvania State University in 2003 and an M.S. and Ph.D. in Aeronautical Engineering from the Air Force Institute of Technology in 2007 and 2018 respectively.  He is a professor with the Department of Aeronautics and Astronautics at the Air Force Institute of Technology.  His research interests include optimal control for vehicle trajectory analysis and dynamics and control of aircraft systems.
\end{biographywithpic}

\begin{biographywithpic}
{Isaac E. Weintraub}{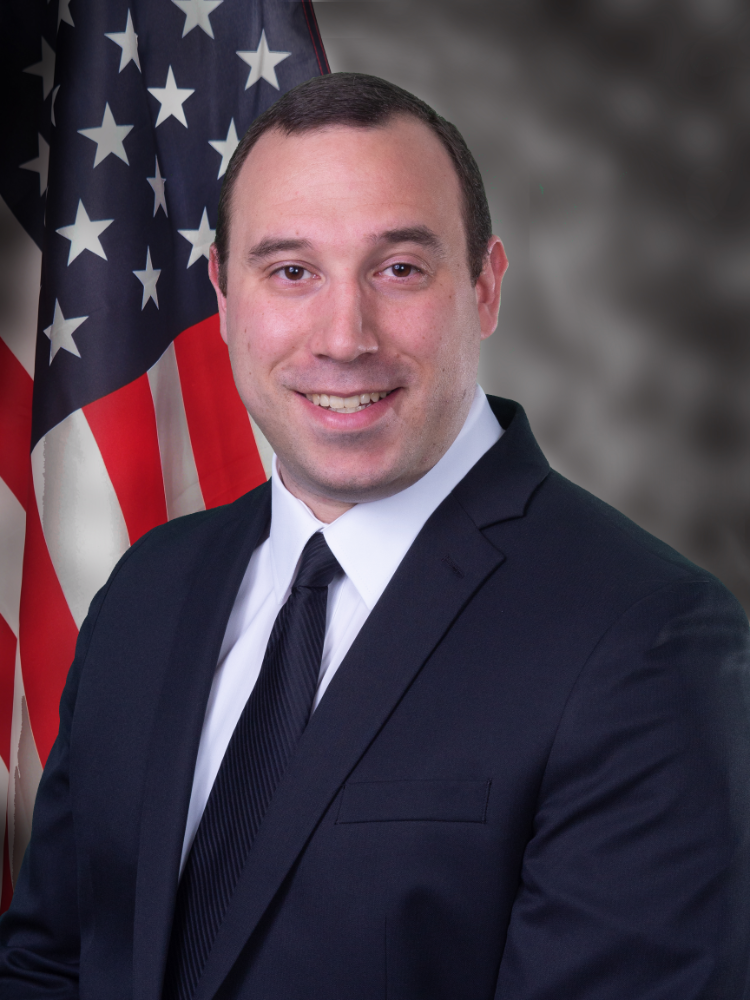}
(S'09-M'15-SM'21) holds a Ph.D. from The Air Force Institute of Technology (2021), an M.S. in Electrical Engineering from University of Texas at Arlington (2011), and a B.S. in Mechanical Engineering from Rose-Hulman Institute of Technology (2009). He is currently an Electronics Engineer with the Control Science Center, Air Force Research Laboratory, Wright-Patterson Air Force Base, Dayton, OH, USA. His research interests are in automation and control of aerospace systems for defense applications.
\end{biographywithpic}

\end{document}